\documentclass[11pt,oneside,reqno]{amsart}
\usepackage[utf8]{inputenc}
\usepackage{amsmath}
\usepackage{amsfonts}
\usepackage{amssymb}
\usepackage{setspace}
\usepackage{amssymb}
\usepackage[shortlabels]{enumitem}
\usepackage{xcolor}
\usepackage{graphicx}
\usepackage{amsthm}
\newtheorem{thm}{Theorem}[section]

\numberwithin{equation}{section}
\usepackage{bbm}
\newtheorem{rmk}{Remark}[section]

\newtheorem{prop}{Proposition}[section]
\newtheorem{lm}{Lemma}[section]
\newtheorem{cor}{Corollary}[section]
\usepackage{geometry}
\usepackage{fancyhdr}
\usepackage{hyperref}
\hypersetup{
	colorlinks=true,
	linkcolor=blue,
	filecolor=magenta,      
	urlcolor=blue,
	pdftitle={Overleaf Example},
	citecolor=red,
	pdfpagemode=FullScreen,
}
\date{}
\begin{document}
	\title[Controllability for stochastic heat equations]{Controllability  for forward stochastic heat equations with dynamic boundary conditions without extra forces}
	\author{Said Boulite\,$^1$, Abdellatif Elgrou\,$^2$, Lahcen Maniar\,$^2$, and Omar Oukdach\,$^3$}
	\begin{abstract}
		In this paper, we continue the study of some controllability issues for the forward stochastic heat equation with dynamic boundary conditions. The main novelty in the present paper consists of considering only one control without extra forces in the noise parts. Under a strong measurability condition, and using a spectral inequality,  we first establish an appropriate observability inequality for the corresponding adjoint system. Then, by the classical duality approach, the null and approximate controllability results are established. 
	\end{abstract}
	
	\maketitle
	
	\footnotetext[1]{Cadi Ayyad University, National School of Applied Sciences, LMDP, UMMISCO (IRD-UPMC), Marrakesh P.B. 575, Morocco. E-mail: \href{s.boulite@uca.ma}{\texttt{s.boulite@uca.ma}},}
	\footnotetext[2]{Cadi Ayyad University, Faculty of Sciences Semlalia, LMDP, UMMISCO (IRD-UPMC), P.B. 2390, Marrakesh, Morocco. E-mail: \href{abdoelgrou@gmail.com}{\texttt{abdoelgrou@gmail.com}}, \href{maniar@uca.ma}{\texttt{maniar@uca.ma}}}
	\footnotetext[3]{\author{Moulay Ismaïl University of Meknes, FST Errachidia,  MSISI Laboratory, AM2CSI Group, BP
			509, Boutalamine, Errachidia, Morocco}. E-mail: \href{omar.oukdach@gmail.com}{\texttt{omar.oukdach@gmail.com}}}
	
	\section{Introduction and main results}
	Let $T>0$, $G\subset\mathbb{R}^N$ ($N\geq2$) be a nonempty open bounded domain with a smooth boundary $\Gamma=\partial G$, and $G_0$ be a nonempty open subset of $G$. $E$ is a measurable subset in $(0,T)$ with a positive Lebesgue measure (i.e., $\mathbf{m}(E)>0$), $\chi_E$ (resp., $\chi_{G_0}$) is the characteristic function of $E$ (resp., $G_0$). For any $s\in[0,T)$, we denote by 
	$$Q_s=(s,T)\times G\,\quad \text{and}\,\quad\Sigma_s=(s,T)\times\Gamma.$$ 
	
	Let $(\Omega,\mathcal{F},\{\mathcal{F}_t\}_{t\geq0},\mathbb{P})$ be a fixed complete filtered probability space on which a one-dimensional standard Brownian motion $W(\cdot)$ is defined such that $\{\mathcal{F}_t\}_{t\geq0}$ is the natural filtration generated by $W(\cdot)$ and augmented by all the $\mathbb{P}$-null sets in $\mathcal{F}$.
	
	For any given Banach space $\mathcal{X}$, the space $C([0,T];\mathcal{X})$ denotes the Banach space of all $\mathcal{X}$-valued continuous functions defined on $[0,T]$. For a sub-sigma algebra $\mathcal{G}\subset\mathcal{F}$, we denote by $L^2_{\mathcal{G}}(\Omega;\mathcal{X})$ the Banach space of all $\mathcal{X}$-valued $\mathcal{G}$-measurable random variables $X$ equipped with the canonical norm
	$$|X|_{L^2_\mathcal{G}(\Omega;\mathcal{X})}=\big(\mathbb{E}\vert X\vert_\mathcal{X}^2\big)^{1/2}.$$
	For any $p,q\in[1,\infty)$, we define $L^p_\mathcal{F}(0,T;L^q(\Omega;\mathcal{X}))$ as the Banach space of all $\mathcal{X}$-valued $\{\mathcal{F}_t\}_{t\geq0}$-adapted stochastic processes $\phi(\cdot)$ endowed with the usual norm 
	$$|\phi|_{L^p_\mathcal{F}(0,T;L^q(\Omega;\mathcal{X}))}=\Bigg(\int_0^T \big(\mathbb{E}|\phi(t)|_\mathcal{X}^q\big)^{p/q}dt\Bigg)^{1/p}.$$
	The special space $L^p_\mathcal{F}(0,T;L^p(\Omega;\mathcal{X}))$ is simply denoted by $L^p_\mathcal{F}(0,T;\mathcal{X})$. The space $L^\infty_\mathcal{F}(0,T;\mathcal{X})$ indicates the Banach space of all $\mathcal{X}$-valued $\{\mathcal{F}_t\}_{t\geq0}$-adapted essentially bounded stochastic processes $\phi(\cdot)$ endowed with the  norm
	\begin{align*}
		|\phi|_{L^\infty_\mathcal{F}(0,T;\mathcal{X})}&=\displaystyle\textnormal{esssup}_{t\in[0,T],\,\,\omega\in\Omega}|\phi(t)|_\mathcal{X}\\
		&=\inf\big\{\mathcal{C}\geq0,\;|\phi(t)|_\mathcal{X}\leq \mathcal{C},\;\textnormal{for a.e.,}\; t\in(0,T),\;\,\mathbb{P}\textnormal{-a.s.}\big\},
	\end{align*}
	and when $\mathcal{X}=\mathbb{R}$, we simply write this space as $L^\infty_\mathcal{F}(0,T)$. Similarly, we define the space $L^\infty_\mathcal{F}(0,T;L^q(\Omega;\mathcal{X}))$. Finally, we denote by $L^2_\mathcal{F}(\Omega;C([0,T];\mathcal{X}))$ the Banach space consisting of all $\mathcal{X}$-valued $\{\mathcal{F}_t\}_{t\geq0}$-adapted continuous stochastic processes $\phi(\cdot)$ equipped with the norm 
	$$|\phi|_{L^2_\mathcal{F}(\Omega;C([0,T];\mathcal{X}))}=\Bigg(\mathbb{E}\bigg[\max_{t\in[0,T]}\vert \phi(t)\vert_\mathcal{X}^2\bigg]\Bigg)^{1/2}.$$
	\noindent
	To study problems with dynamic boundary conditions, we consider the following space
	$$ \mathbb{L}^2 := L^2(G,\,dx)\times L^2(\Gamma,\,d\sigma).$$
	Endowed with the  inner product
	$$ \big\langle (y,y_\Gamma),(z,z_\Gamma)\big\rangle_{\mathbb{L}^2}= \langle y,z\rangle_{L^2(G)} + \langle y_\Gamma,z_\Gamma\rangle_{L^2(\Gamma)} = \int_{G}yz \,dx +  \int_{\Gamma}y_\Gamma z_\Gamma \,d\sigma,$$
	$\mathbb{L}^2$ is  a Hilbert space. Here, $dx$ (resp., $d\sigma$) denoted the Lebesgue (resp., surface) measure in $G$ (resp., on $\Gamma$).\\
	
	In this paper, we study the null and approximate controllability issues of the following forward stochastic heat equation with dynamic boundary conditions
	\begin{equation}\label{1.1}
		\begin{cases}
			\begin{array}{ll}
				dy - \Delta y \,dt = \big[a(t)y+\chi_E(t)\chi_{G_0}(x)u\big] \,dt + b(t)y\,dW(t) &\textnormal{in}\,\,Q_0,\\
				dy_\Gamma-\Delta_\Gamma y_\Gamma \,dt+\partial_\nu y \,dt = a(t)y_\Gamma \,dt+b(t)y_\Gamma \,dW(t) &\textnormal{on}\,\,\Sigma_0,\\
				y_\Gamma(t,x)=y\vert_\Gamma(t,x) &\textnormal{on}\,\,\Sigma_0,\\
				(y,y_\Gamma)\vert_{t=0}=(y_0,y_{\Gamma,0}) &\textnormal{in}\,\,G\times\Gamma,
			\end{array}
		\end{cases}
	\end{equation}
	where $(y_0,y_{\Gamma,0})\in L^2_{\mathcal{F}_0}(\Omega;\mathbb{L}^2)$ is the initial state, $a,b\in L^\infty_{\mathcal{F}}(0,T)$, $u\in L^\infty_{\mathcal{F}}(0,T;L^2(\Omega;L^2(G)))$ is the control function acting only in the small space region $G_0$ and in the time subset $E$. Throughout this paper, we denote by $C$ a positive constant that may change from one place to another. In \eqref{1.1}, $\partial_\nu y = (\nabla y \,\cdot\, \nu)\vert_{\Sigma_0}$ indicates the normal derivative where $\nu$ is the outer unit normal vector at the boundary $\Gamma$. This normal derivative is the coupling term between the bulk and surface equations. For more details and a description of the dynamical boundary conditions (sometimes also called Wentzell boundary conditions) and their physical interpretation, we refer to \cite{cocang2008,Gal15}. Now, we recall the definition of the differential operators on $\Gamma$. The tangential gradient of a function $y_\Gamma:\Gamma\rightarrow\mathbb{R}$  is defined by
	$$
	\nabla_{\Gamma} y_\Gamma(x)=(\partial_{i,\Gamma}y_\Gamma^i)_{i=1}^N=\nabla y(x)- (\partial_{\nu}y)\nu(x)=\left(Id-\nu(x) \nu(x)^{\top}\right) \nabla y(x), \quad x\in \Gamma,
	$$
	where $y$ is an extension of $y_\Gamma$ up to an open neighborhood of $\Gamma$ and $\partial_{i,\Gamma}y_\Gamma^i$ ($i=1,\cdot\cdot\cdot, N$) are the partial tangential derivative of $y_\Gamma$ on $\Gamma$. Note that $\nabla_{\Gamma} y_\Gamma(x)$ is the projection of $\nabla y(x)$ onto the tangent plane $T_{x}\Gamma=\{\xi\in\mathbb{R}^N,\; \xi\cdot\nu(x)=0\}$ at the point $x\in\Gamma$. In particular, for any $x\in \Gamma$, we have $\nabla_{\Gamma}y_\Gamma(x)\cdot\nu(x)= 0$. The  tangential divergence of a function $y_\Gamma:\Gamma\rightarrow\mathbb{R}^N$ is defined by
	$$
	\operatorname{div}_{\Gamma} y_\Gamma=\sum_{i=1}^{N} \partial_{i, \Gamma} y_\Gamma^i, 
	$$
	and the Laplace-Beltrami operator is given by
	$$
	\Delta_{\Gamma} y_\Gamma= \operatorname{div}_{\Gamma}(\nabla_{\Gamma}y_\Gamma).
	$$
	Recall that $H^1(\Gamma)$ and $H^2(\Gamma)$ are also real Hilbert spaces endowed with the following respective norms
	$$ |y_\Gamma|_{H^1(\Gamma)}= \langle y_\Gamma,y_\Gamma \rangle^{\frac{	1}{2}}_{H^1(\Gamma)}, \,\text{with} \,\,\langle y_\Gamma,z_\Gamma \rangle_{H^1(\Gamma)}= \int_{\Gamma}y_\Gamma z_\Gamma \,d\sigma + \int_{\Gamma}\nabla_{\Gamma} y_\Gamma\cdot \nabla_{\Gamma}z_\Gamma \,d\sigma, \quad $$
	and 
	$$ |y_\Gamma|_{H^2(\Gamma)}= \langle y_\Gamma,y_\Gamma \rangle^{\frac{	1}{2}}_{H^2(\Gamma)},\, \text{with}\,\, \langle y_\Gamma,z_\Gamma \rangle_{H^2(\Gamma)}= \int_{\Gamma}y_\Gamma z_\Gamma \,d\sigma +\int_{\Gamma}\Delta_{\Gamma} y_\Gamma\, \Delta_{\Gamma}z_\Gamma \,d\sigma.$$              
	In the subsequent sections, we commonly use the following surface divergence formula
	$$
	\int_\Gamma \triangle_\Gamma y_\Gamma\, z_\Gamma\, \,d\sigma = -\int_\Gamma \nabla_\Gamma y_\Gamma\cdot\nabla_\Gamma z_\Gamma\, \,d\sigma\,,\quad y_\Gamma\in H^2(\Gamma),\,\; z_\Gamma\in H^1(\Gamma).
	$$   
	We also consider the spaces
	$$\mathbb{H}^k=\big\{(y,y_\Gamma)\in H^k(G)\times H^k(\Gamma): \,\, y_\Gamma=y|_\Gamma\big\},\quad k=1,2,$$
	viewed as a subspace of $H^k(G)\times H^k(\Gamma)$ with the natural topology inherited by $H^k(G)\times H^k(\Gamma)$.
	
	The null and approximate controllability issues addressed in this paper  are formulated as follows:
	\begin{enumerate}[1.]
		\item System \eqref{1.1} is said to be null controllable at time $T$ if for any initial state $(y_0,y_{\Gamma,0})\in L^2_{\mathcal{F}_0}(\Omega;\mathbb{L}^2)$, there exists a control $u\in L^\infty_{\mathcal{F}}(0,T;L^2(\Omega;L^2(G)))$ such that the corresponding solution $(y,y_\Gamma)$ satisfies that 
		$$(y(T,\cdot),y_\Gamma(T,\cdot))=(0,0)\;\; \textnormal{in}\,\; G\times\Gamma,\,\,\mathbb{P}\textnormal{-a.s.} 
		$$
		\item System \eqref{1.1} is said to be approximately controllable at time $T$ if for any initial state $(y_0,y_{\Gamma,0})\in L^2_{\mathcal{F}_0}(\Omega;\mathbb{L}^2)$, any final state $(y_T,y_{\Gamma,T})\in L^2_{\mathcal{F}_T}(\Omega;\mathbb{L}^2)$ and for all $\varepsilon>0$, there exists a control $u\in L^\infty_{\mathcal{F}}(0,T;L^2(\Omega;L^2(G)))$ such that the associated solution $(y,y_\Gamma)$ satisfies that
		$$\big|(y(T,\cdot),y_\Gamma(T,\cdot))-(y_T,y_{\Gamma,T})\big|_{L^2_{\mathcal{F}_T}(\Omega;\mathbb{L}^2)}\leq\varepsilon.$$
	\end{enumerate}
	
	The first main result of this paper is stated as follows.
	\begin{thm}\label{thmm1.1}
		For any $T>0$, any nonempty open subset $G_0$ of $G$ and for all $E\subset(0,T)$ such that $\textbf{m}(E)>0$, the system \eqref{1.1} is null controllable at time $T$ and the control $u$ can be chosen so that
		$$|u|^2_{L^\infty_{\mathcal{F}}(0,T;L^2(\Omega;L^2(G)))}\leq C|(y_0,y_{\Gamma,0})|^2_{L^2_{\mathcal{F}_0}(\Omega;\mathbb{L}^2)},$$
		where $C$ is a positive constant depending only on $G$, $G_0$, $E$, $T$, $|a|_{L^\infty_\mathcal{F}(0,T)}$ and $|b|_{L^\infty_\mathcal{F}(0,T)}$.
	\end{thm}
	Unlike null controllability, the following approximate controllability result is proved under a stronger assumption on $E$. 
	\begin{thm}\label{thmm1.2}
		The system \eqref{1.1} is approximately controllable at time $T$ if and only if
		\begin{align}\label{sufneccon}
			\textbf{m}((s,T)\cap E)>0,\;\;\,\textnormal{for all}\,\,s\in[0,T).
		\end{align}
	\end{thm}
	As well known, the controllability problem of \eqref{1.1} may be reduced to an appropriate observability property for the corresponding adjoint system. For this, let us consider the following backward stochastic heat equation
	\begin{equation}\label{backadj}
		\begin{cases}dz+\Delta z \,dt=\big[-a(t) z-b(t) Z\big] \,dt+Z \,dW(t) & \text { in }Q_s, \\ dz_{\Gamma}+\Delta_\Gamma z_\Gamma \,dt-\partial_\nu z \,dt=\big[-a(t) z_{\Gamma}-b(t) \widehat{Z}\big] \,dt+\widehat{Z} \,dW(t) & \text { on }\Sigma_s, \\ z_{\Gamma}(t, x)=z|_\Gamma(t, x) & \text { on }\Sigma_s, \\ (z, z_\Gamma)|_{t=T}=(z_T, z_{\Gamma, T}) & \text { in } G \times \Gamma,\end{cases}
	\end{equation}
	where $(z_T, z_{\Gamma, T}) \in L_{\mathcal{F}_T}^2\left(\Omega ; \mathbb{L}^2\right)$ is the terminal state and $(z, z_{\Gamma}, Z, \widehat{Z})$ is the state variable. 
	
	The following observability inequality is the key tool to establish our null and approximate controllability results.
	\begin{thm}\label{thmm01.3}
		If $\textbf{m}((s,T)\cap E)>0$, for all $s\in[0,T)$, then there exists a positive constant $C$ depending on $G$, $G_0$, $E$, $s$, $T$, $a$ and $b$ such that for any $(z_{T},z_{\Gamma,T})\in L^2_{\mathcal{F}_{T}}(\Omega;\mathbb{L}^2)$, the associated solution $(z,z_\Gamma,Z,\widehat{Z})$ of \eqref{backadj} satisfies that
		\begin{align}\label{1.301}
			\mathbb{E}\vert z(s)\vert_{L^2(G)}^2+\mathbb{E}\vert z_{\Gamma}(s)\vert^2_{L^2(\Gamma)}\leq C\vert\chi_E(t)\chi_{G_0}(x) z\vert^2_{L^1_{\mathcal{F}}(s,T;L^2(\Omega; L^2(G)))}.
		\end{align}
	\end{thm}
	Some remarks are in order.
	\begin{rmk}
		\begin{enumerate}[1.]
			\item In this study,  the null controllability of \eqref{1.1} is proved without any additional assumption on the time control region $E$. However, \eqref{sufneccon} is a sufficient and necessary condition for the approximate controllability of \eqref{1.1}. This shows that for stochastic parabolic equations with dynamic boundary conditions $($as remarked in \cite{lu2011some} for classical boundary conditions$)$, the null controllability does not imply approximate controllability. This is one of the novel phenomena that distinguishes the study of the controllability properties of deterministic and stochastic parabolic equations.
			\item Noting that we have used the same coefficients $a$ and $b$ on the bulk and surface equations in \eqref{1.1}. This is needed in Section \ref{sec3} when we use some spectral properties to derive the ordinary backward stochastic differential equation \eqref{adjsto}.
			\item The controllability of \eqref{1.1} when the potentials $a$ and $b$ depending both on $t$ and $x$ is still an open and challenging problem. In such a case, the observability inequality \eqref{1.301} seems difficult to prove using only one control on the drift part. The authors in \cite{tang2009null} used one more extra control on the diffusion part to handle this issue. 
			\item The null and approximate controllability for general forward stochastic parabolic equations $($when coefficients depend on $t$, $x$ and $\omega$$)$ with dynamic boundary conditions are established in \cite{elgrouDBC} using two extra controls on the diffusion parts. The approach applied to get these controllability results is global Carleman estimates.
		\end{enumerate}
	\end{rmk}
	
	Controllability issues are the subject of much research in the literature. In the context of deterministic parabolic equations with classical boundary conditions, the literature is so large and has been extensively studied in many works, we refer for instance to \cite{lebeau1995controle,zabczyk2020mathematical} and the references therein. For parabolic equations with dynamic boundary conditions, the results are relatively recent and many works have been done in various contexts, we mainly refer to \cite{ACMO20, OuBoMa19,BoMaOuNash,  khoutaibi2020null,KMO22}.
	
	As for stochastic parabolic equations, we find numerous publications concerning the controllability of forward and backward stochastic parabolic equations, see, for instance, \cite{barbu2003carleman,elgrouDBC1d,Preprintelgrou23,elgrouDBC,fadili,liu2014global,lu2011some,lu2021mathematical,san23,yan2018carleman,yansing,zhang2008unique}. In the case of static boundary conditions (Dirichlet, Neumann, and Fourier), the first controllability result for forward stochastic parabolic equations (FSPEs) is obtained in \cite{barbu2003carleman}, where the authors describe a particular set of reachable final states. Later in \cite{tang2009null}, using an extra control on the diffusion part, the authors established controllability results for FSPEs by deriving an appropriate Carleman estimate for backward stochastic parabolic equations. In the same manner, in \cite{yan2018carleman}, by the tool of Carleman estimates, the controllability of FSPEs with Neumann and Fourier boundary conditions is proved. In \cite{elgrouDBC}, by adding two extra controls, the authors studied the controllability of FSPEs with dynamic boundary conditions by developing a new global Carleman estimate for backward stochastic parabolic equations with dynamic boundary conditions. 
	
	We emphasize that the controllability of FSPEs (in the case of static or dynamic boundary conditions) using only one control is still an interesting and challenging question. In \cite{lu2011some}, adapting the Lebeau-Robbiano strategy, the null and approximate controllability in the case of static boundary conditions are established. To the best of our knowledge, the present paper is the first to consider the one-force controllability issue in the context of stochastic parabolic equations with dynamic boundary conditions. This can be done thanks to the spectral inequality proved in \cite{spinq}.\\
	
	The paper is organized as follows: In Section \ref{sec2}, we briefly present some preliminaries needed in the sequel. In Section \ref{sec3}, we use a spectral inequality to prove the needed observability inequality \eqref{1.301}. Section \ref{sec4} is devoted to proving the announced controllability results. In Section \ref{sec5}, we end the paper with a conclusion where we discuss some open problems.
	\section{Well-posedness and  preliminary results}\label{sec2}
	This section presents well-posedness results for solutions of the forward and backward stochastic parabolic equations \eqref{1.1} and \eqref{backadj}. Additionally, we provide the spectrum and an important spectral inequality for the operator governing these systems. For more details on well-posedness and the regularity of solutions of forward and backward stochastic partial differential equations, we refer the readers to \cite{da2014stochastic,krylov,lu2021mathematical,pardoux1990adapted}.
	
	Introduce the following linear operator $\mathcal{A}$ associated to the  systems \eqref{1.1} and \eqref{backadj}
	$$\mathcal{A}=\begin{pmatrix} 
		\Delta & 0 \\
		-\partial_{\nu} & \Delta_\Gamma
	\end{pmatrix},\quad\textnormal{with domain}\quad\mathcal{D}(\mathcal{A})=\mathbb{H}^2.$$
	We have the following generation result. For the proof, see \cite{maniar2017null}.  
	\begin{prop}\label{prop2.1}
		The operator $\mathcal{A}$ is a densely defined, self-adjoint, dissipative operator and generates an analytic $C_0$-semigroup $(e^{t\mathcal{A}})_{t\geq0}$ on $\mathbb{L}^2$.
	\end{prop}
	\noindent
	Observe that the system \eqref{1.1} can be written in the following abstract Cauchy problem
	\begin{equation}\label{abs1.1}
		\begin{cases}
			d\mathcal{Y}=(\mathcal{A}\mathcal{Y}+B_1(t)\mathcal{Y} + F_1(t)) \,dt + B_2(t)\mathcal{Y} \,dW(t)\,,\\
			\mathcal{Y}(0)=\mathcal{Y}_0\in L^2_{\mathcal{F}_0}(\Omega;\mathbb{L}^2),
	\end{cases}\end{equation}
	where 
	$$\mathcal{Y}=
	\begin{pmatrix}
		y\\
		y_\Gamma
	\end{pmatrix},\; \mathcal{Y}_0=
	\begin{pmatrix}
		y_0\\
		y_{\Gamma,0}
	\end{pmatrix},\; B_1=
	\begin{pmatrix}
		a&0\\
		0&a
	\end{pmatrix},\; B_2=
	\begin{pmatrix}
		b&0\\
		0&b
	\end{pmatrix},\; F_1=
	\begin{pmatrix}
		\chi_E\chi_{G_0}u\\
		0
	\end{pmatrix}.
	$$
	By using Proposition \ref{prop2.1} and Theorem 3.24 in \cite{lu2021mathematical}, we deduce the following well-posedness for the system \eqref{1.1}.
	\begin{thm}\label{thm2.1}
		For each  $(y_0,y_{\Gamma,0})\in L^2_{\mathcal{F}_0}(\Omega;\mathbb{L}^2)$ and $u\in L^\infty_{\mathcal{F}}(0,T;L^2(\Omega;L^2(G)))$, the system \eqref{abs1.1} $($and hence \eqref{1.1}$)$ admits a unique mild solution
		$$\mathcal{Y}=(y,y_\Gamma)\in L^2_\mathcal{F}(\Omega;C([0,T];\mathbb{L}^2))\bigcap L^2_\mathcal{F}(0,T;\mathbb{H}^1)$$
		given by
		$$\mathcal{Y}(t)=e^{t\mathcal{A}}\mathcal{Y}_0+\int_0^t e^{(t-\tau)\mathcal{A}}\big[B_1(\tau)\mathcal{Y}(\tau)+F_1(\tau)\big]\,d\tau+\int_0^t e^{(t-\tau)\mathcal{A}} B_2(\tau) \,dW(\tau)$$
		for all $t\in[0,T]$, \,$\mathbb{P}$\textnormal{-a.s.} Moreover,
		\begin{align*}
			\vert(y,y_\Gamma)\vert_{L^2_\mathcal{F}(\Omega;C([0,T];\mathbb{L}^2))} + \vert(y,y_\Gamma)\vert_{L^2_\mathcal{F}(0,T;\mathbb{H}^1)}\leq C\,\big(|(y_0,y_{\Gamma,0})|_{L^2_{\mathcal{F}_0}(\Omega;\mathbb{L}^2)}+|u|_{L^\infty_{\mathcal{F}}(0,T;L^2(\Omega;L^2(G)))}\big).
		\end{align*}
	\end{thm}
	\noindent
	In the same way,  we can rewrite equation \eqref{backadj} in the following abstract Cauchy problem
	\begin{equation}\label{absback}
		\begin{cases}
			d\mathcal{Z}=(-\mathcal{A}\mathcal{Z}+B_3(t)\mathcal{Z}+B_4(t)\widehat{\mathcal{Z}}) \,dt + \widehat{\mathcal{Z}} \,dW(t)\,,\\
			\mathcal{Z}(T)=\mathcal{Z}_T\in L^2_{\mathcal{F}_T}(\Omega;\mathbb{L}^2),
	\end{cases}\end{equation}
	where $$\mathcal{Z}=
	\begin{pmatrix}
		z\\
		z_\Gamma
	\end{pmatrix},\; \mathcal{Z}_T=
	\begin{pmatrix}
		z_T\\
		z_{\Gamma,T}
	\end{pmatrix},\;
	\widehat{\mathcal{Z}}=
	\begin{pmatrix}
		Z\\
		\widehat{Z}
	\end{pmatrix},\; B_3(t)=
	\begin{pmatrix}
		-a&0\\
		0&-a
	\end{pmatrix},\; B_4(t)=
	\begin{pmatrix}
		-b&0\\
		0&-b
	\end{pmatrix}.
	$$
	From Proposition \ref{prop2.1} and Theorem 4.11 in \cite{lu2021mathematical}, we have the following well-posedness result for the system \eqref{backadj}.
	\begin{thm}
		For each  $(z_T, z_{\Gamma,T}) \in L_{\mathcal{F}_T}^2\left(\Omega ; \mathbb{L}^2\right)$, the system \eqref{absback} $($and hence \eqref{backadj}$)$ has a unique mild solution
		$$
		(\mathcal{Z},\widehat{\mathcal{Z}})=(z, z_{\Gamma}, Z, \widehat{Z}) \in\left(L_{\mathcal{F}}^2\left(\Omega ; C\left([s, T] ; \mathbb{L}^2\right)\right) \bigcap L_{\mathcal{F}}^2\left(s, T ; \mathbb{H}^1\right)\right) \times L_{\mathcal{F}}^2\left(s, T ; \mathbb{L}^2\right)
		$$
		given by
		$$\mathcal{Z}(t)=e^{(T-t)\mathcal{A}}\mathcal{Z}_T-\int_t^T e^{(\tau-t)\mathcal{A}} \big[B_3(\tau)\mathcal{Z}(\tau)+B_4(\tau)\widehat{\mathcal{Z}}(\tau)\big]\,d\tau-\int_t^T e^{(\tau-t)\mathcal{A}} \widehat{\mathcal{Z}}(\tau) \,dW(\tau)$$
		for all $t\in[s,T]$, \,$\mathbb{P}$\textnormal{-a.s.} Furthermore, it holds that
		\begin{align}\label{wbac2.1}
			\left|\left(z, z_{\Gamma}\right)\right|_{L_{\mathcal{F}}^2\left(\Omega ; C\left([s, T] ; \mathbb{L}^2\right)\right)}+\left|\left(z, z_{\Gamma}\right)\right|_{L_{\mathcal{F}}^2\left(s, T ; \mathbb{H}^1\right)}+|(Z, \widehat{Z})|_{L_{\mathcal{F}}^2\left(s, T ; \mathbb{L}^2\right)} \leq C\left|\left(z_T, z_{\Gamma, T}\right)\right|_{L_{\mathcal{F}_T}^2\left(\Omega ; \mathbb{L}^2\right)}.
		\end{align}
	\end{thm}
	We end this section with two important spectral results that will be useful for the rest of the paper. Firstly, we give the form of the spectrum of $\mathcal{A}$. For the proof, see \cite{Gal15}.
	\begin{lm}\label{lm2.2}
		There exists a sequence of numbers $$ 0\leq\lambda_1\leq \lambda_2\leq \cdots\leq\lambda_k\leq\lambda_{k+1}\leq \cdots, $$
		converging to $+\infty$ such that the spectrum of $-\mathcal{A}$ is given by 
		$$\sigma(-\mathcal{A})=\{\lambda_k:\quad k\geq1 \},$$ with $\Psi_k = (\psi_k , \psi_{\Gamma,k}) \in \mathbb{H}^2$ $($for $k=1,2,\cdot\cdot\cdot$$)$ be the corresponding eigenfunctions which constitute a Hilbert basis of the space $\mathbb{L}^2$.
	\end{lm}
	Secondly, we recall the following spectral inequality. For the proof, we refer to \cite{spinq}. 
	\begin{lm}\label{Spinq}
		There exists a constant $C>0$ such that
		\begin{equation}\label{specine}
			\sqrt{\sum_{\lambda_{j} \leq r}\left|a_{j}\right|^{2}} \leq C e^{C \sqrt{r}}\Bigg| \sum_{\lambda_{j} \leq r} a_{j} \psi_{j}\Bigg|_{L^{2}(G_0)},
		\end{equation}
		for any $r>0$ and all sequence $\left(a_{j}\right)_{j \geq 1} \subset \mathbb{R}$.
	\end{lm}
	\begin{rmk}
		The above spectral inequality implies the following partial observability inequality: There exists $C>0$ such that for all $r>0$, we have 
		\begin{equation}\label{sp2}
			| z|^2_{L^2(G)} + | z_{\Gamma }|^2_{L^2(\Gamma)}\leq C e^{C \sqrt{r}}\int_{G_0} z^2 \,dx,   \quad \forall (z, z_{\Gamma})\in V_{r},
		\end{equation} 
		where $V_{r}=\textnormal{span}\{\Psi_j\in \mathbb{L}^2 : \;\, \lambda_j\leq r\}$. This inequality can show a partial null controllability result of \eqref{1.1}. That is, for all initial state $(y_0,y_{\Gamma,0})\in L^2_{\mathcal{F}_0}(\Omega;\mathbb{L}^2)$, there exists a control $u\in L^2_\mathcal{F}(0,T; L^2(G))$ so that the associated solution $(y,y_\Gamma)$ of \eqref{1.1} satisfies
		\begin{equation*}
			\begin{cases}
				\Pi_{V_r}(y(T,\cdot),y_\Gamma(T,\cdot))=(0,0)\;\textnormal{in}\;G\times\Gamma,\;\mathbb{P}\textnormal{-a.s.},\\
				|u|^2_{L^2_\mathcal{F}(0,T;L^2(G))}\leq \frac{C}{T}e^{C\sqrt{r}}|(y_0,y_{\Gamma,0})|^2_{L^2_{\mathcal{F}_0}(\Omega;\mathbb{L}^2)},
			\end{cases}
		\end{equation*}
		where $\Pi_{V_r}$ is the orthogonal projection from $\mathbb{L}^2$ to $V_r$.
	\end{rmk}
	\section{Observability problem}\label{sec3}
	This section is devoted to proving the observability inequality \eqref{1.301} for the adjoint backward stochastic heat equation \eqref{backadj}. In what follows, we set
	\begin{align}\label{3.2delta}
		\delta=2\vert a\vert_{L^\infty_{\mathcal{F}}(0,T)}+\vert b\vert^2_{L^\infty_{\mathcal{F}}(0,T)}.
	\end{align}
	Let $r>0$ and for any $j\geq1$,  $\lambda_j$ and  $\Psi_j= (\psi_j,\psi_{\Gamma,j})$ are the eigen-elements of the operator $-\mathcal{A}$. For any $(f,g)\in\mathbb{L}^2$,  we have that
	$$(f,g)=\mathcal{E}_r(f,g) + \mathcal{E}^\bot_r(f,g),$$
	where 
	$$\mathcal{E}_r(f,g)=\sum_{\lambda_j\leq r}\langle(f,g),(\psi_j,\psi_{\Gamma,j})\rangle_{\mathbb{L}^2}(\psi_j,\psi_{\Gamma,j}),
	$$
	and
	$$\mathcal{E}^\bot_r(f,g)=\sum_{\lambda_j>r}\langle(f,g),(\psi_j,\psi_{\Gamma,j})\rangle_{\mathbb{L}^2}(\psi_j,\psi_{\Gamma,j}).$$
	Let $(z,z_\Gamma,Z,\widehat{Z})[T,\cdot,(z_T,z_{\Gamma,T})]$ be the solution of \eqref{backadj} with the terminal state $(z_T,z_{\Gamma,T})$ at $t=T$ and $z^T_j=\langle(z_T,z_{\Gamma,T}),(\psi_j,\psi_{\Gamma,j})\rangle_{\mathbb{L}^2}$. Then, it is easy to see that for any $t\in(s,T)$
	\begin{align}\label{zz1for}
		\mathcal{E}_r(z(t),z_\Gamma(t))=\sum_{\lambda_j\leq r}z_j(t,T,z^T_j)(\psi_j,\psi_{\Gamma,j})=(z,z_\Gamma)[t,T,\mathcal{
			E}_r(z_T,z_{\Gamma,T})],
	\end{align}
	where $(z,z_\Gamma,Z,\widehat{Z})$ is the solution of \eqref{backadj} with the terminal state $\mathcal{E}_r(z_T,z_{\Gamma,T})$ and
	\begin{align}\label{zz2for}
		\mathcal{E}^\bot_r(z(t),z_\Gamma(t))=\sum_{\lambda_j>r}z_j(t,T,z^T_j)(\psi_j,\psi_{\Gamma,j})=(z,z_\Gamma)[t,T,\mathcal{
			E}_r^\bot(z_T,z_{\Gamma,T})],
	\end{align}
	where $(z,z_\Gamma,Z,\widehat{Z})$ is the solution of \eqref{backadj} with the terminal state $\mathcal{E}_r^\bot(z_T,z_{\Gamma,T})$.
	
	In \eqref{zz1for} and \eqref{zz2for}, the pair $(z_j,Z_j)[\cdot,T,z^T_j]$ is the solution of the following backward stochastic differential equation
	\begin{equation}\label{adjsto}
		\begin{cases}
			d z_j-\lambda_j z_j \,dt = \big[-a(t) z_j - b(t)Z_j\big] \,dt+Z_j \,dW(t) & \text {in } [s, T), \\
			z_j(T)=z^T_j. & 
		\end{cases}
	\end{equation}
	
	Firstly, let us show the following interpolation inequality.
	\begin{prop}\label{propp3.11}
		There exists a constant $C>0$ such that for any $(z_T,z_{\Gamma,T})\in L^2_{\mathcal{F}_T}(\Omega;\mathbb{L}^2)$ and  $t\in[s,T)$, we have
		\begin{align}\label{interineq}
			\mathbb{E}|(z(t),z_\Gamma(t))|^2_{\mathbb{L}^2}\leq Ce^{C(T-t)^{-1}}\big(\mathbb{E}|z(t)|^2_{L^2(G_0)}\big)^{1/2}\big(\mathbb{E}|(z_T,z_{\Gamma,T})|^2_{\mathbb{L}^2}\big)^{1/2},
		\end{align}
		where $(z,z_\Gamma,Z,\widehat{Z})$ is the solution of \eqref{backadj} with final state $(z_T,z_{\Gamma,T})$.
	\end{prop}
	\begin{proof}
		From the spectral inequality \eqref{specine}, it is easy to see that for any $t\in[s,T)$
		\begin{align*}
			|\mathcal{E}_r(z(t),z_\Gamma(t))|^2_{\mathbb{L}^2}\leq Ce^{C\sqrt{r}}\Bigg|\sum_{\lambda_j\leq r}z_j(t)\psi_j\Bigg|^2_{L^2(G_0)},
		\end{align*}
		which implies that
		\begin{align}\label{3.9901}
			|\mathcal{E}_r(z(t),z_\Gamma(t))|^2_{\mathbb{L}^2}\leq Ce^{C\sqrt{r}}\big(|z(t)|_{L^2(G_0)}^2+|\mathcal{E}^\bot_r (z(t),z_\Gamma(t))|^2_{\mathbb{L}^2}\big).
		\end{align}
		Note also that
		\begin{align}\label{3.10012}
			|(z(t),z_\Gamma(t))|^2_{\mathbb{L}^2}=|\mathcal{E}_r(z(t),z_\Gamma(t))|^2_{\mathbb{L}^2}+|\mathcal{E}_r^\bot(z(t),z_\Gamma(t))|^2_{\mathbb{L}^2}.
		\end{align}
		Combining \eqref{3.9901} and \eqref{3.10012}, we deduce that
		\begin{align}\label{firstine}
			\mathbb{E}|(z(t),z_\Gamma(t))|^2_{\mathbb{L}^2}\leq 2Ce^{C\sqrt{r}}\big(\mathbb{E}|z(t)|_{L^2(G_0)}^2+\mathbb{E}|\mathcal{E}^\bot_r (z(t),z_\Gamma(t))|^2_{\mathbb{L}^2}\big).
		\end{align}
		On the other hand, by Itô's formula, we compute $d_\tau(e^{(2r-\delta)(T-\tau)}|\mathcal{E}^\bot_r(z(\tau),z_\Gamma(\tau))|^2_{\mathbb{L}^2})$, integrating the obtained equality w.r.t $\tau\in(t,T)$ and taking the expectation on both sides, we obtain that
		\begin{align}\label{3.4401}
			\begin{aligned}
				&\,\mathbb{E}\vert\mathcal{E}_r^\bot(z_T,z_{\Gamma,T})\vert^2_{\mathbb{L}^2}-e^{(2r-\delta)(T-t)}\mathbb{E}|\mathcal{E}^\bot_r(z(t),z_\Gamma(t))|^2_{\mathbb{L}^2}\\
				&=-(2r-\delta)\mathbb{E}\int_t^T e^{(2r-\delta)(T-\tau)}|(z,z_\Gamma)|_{\mathbb{L}^2}^2\,d\tau\\
				&\quad\,+\mathbb{E}\int_t^T e^{(2r-\delta)(T-\tau)} \bigg[2\sum_{\lambda_j>r}\lambda_j z_j^2-2a(\tau)\vert(z,z_\Gamma)\vert^2_{\mathbb{L}^2}\\
				&\quad\,-2b(\tau)\langle(z,z_\Gamma),(Z,\widehat{Z})\rangle_{\mathbb{L}^2}+|(Z,\widehat{Z})|^2_{\mathbb{L}^2}\bigg]\,d\tau.
			\end{aligned}
		\end{align}
		Using Young inequality for the fourth term on the right-hand side of \eqref{3.4401} and  the fact  that $|(z,z_\Gamma)|_{\mathbb{L}^2}^2=\sum_{\lambda_j>r}z_j^2$, we get  
		\begin{align*}
			\begin{aligned}
				&\,\mathbb{E}\vert\mathcal{E}_r^\bot(z_T,z_{\Gamma,T})\vert^2_{\mathbb{L}^2}-e^{(2r-\delta)(T-t)}\mathbb{E}|\mathcal{E}^\bot_r(z(t),z_\Gamma(t))|^2_{\mathbb{L}^2}\\
				&\geq-(2r-\delta)\mathbb{E}\int_t^T e^{(2r-\delta)(T-\tau)}\sum_{\lambda_j>r}z_j^2\,d\tau\\
				&\quad\,+\mathbb{E}\int_t^T e^{(2r-\delta)(T-\tau)} \bigg[2\sum_{\lambda_j>r}\lambda_j z_j^2-2|a(\tau)|\vert(z,z_\Gamma)\vert^2_{\mathbb{L}^2}-|b(\tau)|^2\vert(z,z_\Gamma)\vert^2_{\mathbb{L}^2}\bigg]\,d\tau.
			\end{aligned}
		\end{align*}
		Then, it follows that
		\begin{align*}
			&\,\mathbb{E}\vert\mathcal{E}_r^\bot(z_T,z_{\Gamma,T})\vert^2_{\mathbb{L}^2}-e^{(2r-\delta)(T-t)}\mathbb{E}|\mathcal{E}^\bot_r(z(t),z_\Gamma(t))|^2_{\mathbb{L}^2}\\
			&\geq \mathbb{E}\int_t^T (\delta-2|a(\tau)|-|b(\tau)|^2) e^{(2r-\delta)(T-\tau)}|(z,z_\Gamma))|^2_{\mathbb{L}^2}\,d\tau.
		\end{align*}
		Recalling \eqref{3.2delta}, we deduce 
		\begin{align*}
			\mathbb{E}\vert\mathcal{E}_r^\bot(z_T,z_{\Gamma,T})\vert^2_{\mathbb{L}^2}-e^{(2r-\delta)(T-t)}\mathbb{E}|\mathcal{E}^\bot_r(z(t),z_\Gamma(t))|^2_{\mathbb{L}^2}\geq0,
		\end{align*}
		which implies  that
		\begin{align}\label{3.33012}
			\mathbb{E}|\mathcal{E}^\bot_r(z(t),z_\Gamma(t))|^2_{\mathbb{L}^2} \leq e^{(-2r+\delta)(T-t)}\mathbb{E}\vert(z_T,z_{\Gamma,T})\vert^2_{\mathbb{L}^2}.
		\end{align}
		Now, combining \eqref{3.33012} and \eqref{firstine}, we conclude that
		\begin{align*}
			\mathbb{E}|(z(t),z_\Gamma(t))|^2_{\mathbb{L}^2}\leq 2Ce^{C\sqrt{r}}(\mathbb{E}|z(t)|_{L^2(G_0)}^2+e^{(-2r+\delta)(T-t)}\mathbb{E}\vert(z_T,z_{\Gamma,T})\vert^2_{\mathbb{L}^2}).
		\end{align*}
		Hence
		\begin{align}\label{3.110125}
			\mathbb{E}|(z(t),z_\Gamma(t))|^2_{\mathbb{L}^2}\leq 2Ce^{\delta T}e^{C\sqrt{r}-r(T-t)}\big[e^{r(T-t)}\mathbb{E}|z(t)|_{L^2(G_0)}^2+e^{-r(T-t)}\mathbb{E}\vert(z_T,z_{\Gamma,T})\vert^2_{\mathbb{L}^2}\big].
		\end{align}
		Using the fact that
		$$C\sqrt{r}-r(T-t)\leq\frac{C^2}{4(T-t)},$$ and the inequality \eqref{3.110125}, we deduce for some $C>0$ that
		\begin{align*}
			\mathbb{E}|(z(t),z_\Gamma(t))|^2_{\mathbb{L}^2}\leq Ce^{C(T-t)^{-1}}\big[e^{r(T-t)}\mathbb{E}|z(t)|_{L^2(G_0)}^2+e^{-r(T-t)}\mathbb{E}\vert(z_T,z_{\Gamma,T})\vert^2_{\mathbb{L}^2}\big].
		\end{align*}
		Then, we get 
		\begin{align*}
			\mathbb{E}|(z(t),z_\Gamma(t))|^2_{\mathbb{L}^2}\leq Ce^{C(T-t)^{-1}}\big[\varepsilon^{-1}\mathbb{E}|z(t)|_{L^2(G_0)}^2+\varepsilon\mathbb{E}\vert(z_T,z_{\Gamma,T})\vert^2_{\mathbb{L}^2}\big],
		\end{align*}
		for any $\varepsilon\in(0,1)$. Therefore, we end up with 
		\begin{align}\label{3.1312012}
			\mathbb{E}|(z(t),z_\Gamma(t))|^2_{\mathbb{L}^2}\leq Ce^{C(T-t)^{-1}}\big[\varepsilon^{-1}A+\varepsilon A^{-1}\big]\times(\mathbb{E}|z(t)|^2_{L^2(G_0)})^{1/2}(\mathbb{E}|(z_T,z_{\Gamma,T})|^2_{\mathbb{L}^2})^{1/2},
		\end{align}
		where $A=\big(\mathbb{E}|z(t)|^2_{L^2(G_0)}\big)^{1/2}\big(\mathbb{E}|(z_T,z_{\Gamma,T})|^2_{\mathbb{L}^2}\big)^{-1/2}.$ Now, from  \eqref{wbac2.1}, it is easy to see that there exists a positive constant $C$ such that $A\leq C$. Finally, it is not difficult to see that there exists a small $\varepsilon\in(0,1)$ and a constant $C>0$ such that $\varepsilon^{-1}A+\varepsilon A^{-1}\leq C$. This implies the desired estimate \eqref{interineq} and then completes the proof of Proposition \ref{propp3.11}.
	\end{proof}
	
	Now, we are in a position to prove the observability inequality \eqref{1.301}. The proof borrows some ideas from \cite{observineqback}.
	\begin{proof}[Proof of Theorem \ref{thmm01.3}]
		Let $\tilde{t}\in(s,T)\cap E$, according to \cite{Lionsopt71}, we have that for any $\rho\in(0,1)$, there exists an increasing sequence $\{t_i\}_{i\geq1}\subset(s,T)$ such that
		\begin{align}\label{p1}
			t_1<t_2<\cdot\cdot\cdot<t_i<t_{i+1}<\cdot\cdot\cdot<\tilde{t}\;\;\textnormal{and}\;\;t_i\rightarrow\tilde{t}\,\,\textnormal{as}\,\,i\rightarrow\infty,
		\end{align}
		\begin{align}\label{p2}
			\textbf{m}(E\cap(t_i,t_{i+1}))\geq \frac{1}{3}(t_{i+1}-t_i),\quad\forall i\geq1,
		\end{align}
		\begin{align}\label{p3}
			t_{i+2}-t_{i+1}=\rho(t_{i+1}-t_i),\quad\forall i\geq1.
		\end{align}
		Let us fix
		\begin{align}\label{p4}
			\eta_i=t_{i+1}-\frac{1}{6}(t_{i+1}-t_i),\quad\forall i\geq1.
		\end{align}
		From the interpolation inequality \eqref{interineq}, we deduce that for any $i\geq1$ and for all $t\in(t_i,\eta_i)$
		\begin{align}\label{3.191901}
			\mathbb{E}|(z(t),z_\Gamma(t))|^2_{\mathbb{L}^2}\leq Ce^{C(t_{i+1}-t)^{-1}}\big(\mathbb{E}|z(t)|^2_{L^2(G_0)}\big)^{1/2}\big(\mathbb{E}|(z(t_{i+1}),z_{\Gamma}(t_{i+1}))|^2_{\mathbb{L}^2}\big)^{1/2}.
		\end{align}
		It is easy to see that there exists a constant $C>0$ so that for any $t\in(t_i,\eta_i)$
		$$\mathbb{E}|(z(t_i),z_\Gamma(t_i))|^2_{\mathbb{L}^2}\leq C\,\mathbb{E}|(z(t),z_\Gamma(t))|^2_{\mathbb{L}^2},$$
		and also
		$$t_{i+1}-t\geq\frac{1}{6}(t_{i+1}-t_i).$$
		Then, it follows from \eqref{3.191901} that
		\begin{align}\label{3.1202001}
			\mathbb{E}|(z(t_i),z_\Gamma(t_i))|^2_{\mathbb{L}^2}\leq Ce^{C(t_{i+1}-t_i)^{-1}}\big(\mathbb{E}|z(t)|^2_{L^2(G_0)}\big)^{1/2}\big(\mathbb{E}|(z(t_{i+1}),z_{\Gamma}(t_{i+1}))|^2_{\mathbb{L}^2}\big)^{1/2}.
		\end{align}
		Applying Young's inequality on the right-hand side of \eqref{3.1202001}, we deduce that for any $\varepsilon>0$ and any $t\in(t_i,\eta_i)$, we have
		\begin{align*}
			\mathbb{E}|(z(t_i),z_\Gamma(t_i))|^2_{\mathbb{L}^2}\leq \varepsilon^{-2}Ce^{C(t_{i+1}-t_i)^{-1}}\mathbb{E}|z(t)|^2_{L^2(G_0)}+\varepsilon^2\mathbb{E}|(z(t_{i+1}),z_{\Gamma}(t_{i+1}))|^2_{\mathbb{L}^2},
		\end{align*}
		which can be written as
		\begin{align}\label{3.121212101}
			\Lambda_i\leq \varepsilon^{-1}Ce^{C(t_{i+1}-t_i)^{-1}}\big(\mathbb{E}|z(t)|^2_{L^2(G_0)}\big)^{1/2}+\varepsilon\Lambda_{i+1},
		\end{align}
		where $\Lambda_i=\big(\mathbb{E}|(z(t_i),z_\Gamma(t_i))|^2_{\mathbb{L}^2}\big)^{1/2}$. Integrating \eqref{3.121212101} on $E\cap(t_i,\eta_i)$, we get that
		\begin{align}\label{3.22220123}
			\textbf{m}(E\cap(t_i,\eta_i))\Lambda_i\leq \varepsilon^{-1}Ce^{C(t_{i+1}-t_i)^{-1}}\int_{t_i}^{t_{i+1}} \chi_E(t)\big(\mathbb{E}|z(t)|^2_{L^2(G_0)}\big)^{1/2}dt+\varepsilon \textbf{m}(E\cap(t_i,\eta_i))\Lambda_{i+1}.
		\end{align}
		On the other hand, from \eqref{p2} and \eqref{p4}, it is easy to check that
		\begin{align}\label{inemes}
			\textbf{m}(E\cap(t_i,\eta_i))=\textbf{m}(E\cap(t_i,t_{i+1}))-\textbf{m}(E\cap(\eta_i,t_{i+1}))\geq\frac{1}{6}(t_{i+1}-t_i).
		\end{align}
		From \eqref{inemes} and \eqref{3.22220123}, we obtain that
		\begin{align*}
			\Lambda_i\leq \varepsilon^{-1}Ce^{C(t_{i+1}-t_i)^{-1}}\int_{t_i}^{t_{i+1}} \chi_E(t)\big(\mathbb{E}|z(t)|^2_{L^2(G_0)}\big)^{1/2}dt+\varepsilon\Lambda_{i+1}.
		\end{align*}
		Then, it follows that
		\begin{align}\label{ine3.2421}
			\varepsilon e^{-C(t_{i+1}-t_i)^{-1}}\Lambda_i\leq \varepsilon^2 e^{-C(t_{i+1}-t_i)^{-1}}\Lambda_{i+1}+ C\int_{t_i}^{t_{i+1}} \chi_E(t)\big(\mathbb{E}|z(t)|^2_{L^2(G_0)}\big)^{1/2}dt.
		\end{align}
		By choosing $\varepsilon=e^{-\frac{1}{2}(t_{i+1}-t_i)^{-1}}$, the inequality \eqref{ine3.2421} implies that
		\begin{align}\label{est3.25250}
			e^{-(C+\frac{1}{2})(t_{i+1}-t_i)^{-1}}\Lambda_i\leq e^{-(C+1)(t_{i+1}-t_i)^{-1}}\Lambda_{i+1}+ C\int_{t_i}^{t_{i+1}} \chi_E(t)\big(\mathbb{E}|z(t)|^2_{L^2(G_0)}\big)^{1/2}dt.
		\end{align}
		Recalling \eqref{p3}, then \eqref{est3.25250} provides that
		\begin{align}\label{est3.262601}
			e^{-(C+\frac{1}{2})(t_{i+1}-t_i)^{-1}}\Lambda_i\leq e^{-(C+1)\rho(t_{i+2}-t_{i+1})^{-1}}\Lambda_{i+1}+ C\int_{t_i}^{t_{i+1}} \chi_E(t)\big(\mathbb{E}|z(t)|^2_{L^2(G_0)}\big)^{1/2}dt.
		\end{align}
		Now, by taking $\rho=(C+1/2)/(C+1)$ in \eqref{est3.262601}, we end up with
		\begin{align}\label{3.27270102}
			e^{-(C+\frac{1}{2})(t_{i+1}-t_i)^{-1}}\Lambda_i\leq e^{-(C+\frac{1}{2})(t_{i+2}-t_{i+1})^{-1}}\Lambda_{i+1}+ C\int_{t_i}^{t_{i+1}} \chi_E(t)\big(\mathbb{E}|z(t)|^2_{L^2(G_0)}\big)^{1/2}dt.
		\end{align}
		Since the estimate \eqref{3.27270102} holds for all $i\geq1$, then by summing all these inequalities from $i=1$ to $+\infty$ and using \eqref{p1}, we conclude that
		\begin{align*}
			e^{-(C+\frac{1}{2})(t_{2}-t_1)^{-1}}\Lambda_1\leq C\int_{t_1}^{\tilde{t}} \chi_E(t)\big(\mathbb{E}|z(t)|^2_{L^2(G_0)}\big)^{1/2}dt.
		\end{align*}
		It follows that
		\begin{align}\label{eqq3.26}
			\mathbb{E}|(z(t_1),z_\Gamma(t_1))|_{\mathbb{L}^2}^2\leq Ce^{(C+1)(t_{2}-t_1)^{-1}}\Bigg(\int_{t_1}^{\tilde{t}} \chi_E(t)\big(\mathbb{E}|z(t)|^2_{L^2(G_0)}\big)^{1/2}dt\Bigg)^2.
		\end{align}
		On the other hand by \eqref{wbac2.1}, it is easy to see that there exists a constant $C>0$ so that 
		\begin{align}\label{eq3.277}
			\mathbb{E}|(z(s),z_\Gamma(s))|^2_{\mathbb{L}^2}\leq C \,\mathbb{E}|(z(t_1),z_\Gamma(t_1))|^2_{\mathbb{L}^2}.
		\end{align}
		Finally, from \eqref{eqq3.26} and \eqref{eq3.277}, we deduce
		$$\mathbb{E}|(z(s),z_\Gamma(s))|_{\mathbb{L}^2}^2 \leq C \Bigg(\int_s^T \big(\mathbb{E}|\chi_E(t)\chi_{G_0}(x) z|^2_{L^2(G)}\big)^{1/2}dt\Bigg)^2,$$
		which immediately implies the desired observability inequality \eqref{1.301}. This concludes the proof of Theorem \ref{thmm01.3}.
	\end{proof}
	Now, using Theorem \ref{thmm01.3}, it is easy to deduce the following unique continuation property for the system \eqref{backadj}.
	\begin{cor}\label{corunconprop}
		If $\textbf{m}((s,T)\cap E)>0$, for any $s\in[0,T)$, then the solutions of \eqref{backadj} satisfy that
		\begin{align*}
			z=0\,\,\;\textnormal{in}\,\,E\times G_0,\,\;\mathbb{P}\textnormal{-a.s.} \,\,\Longrightarrow\,\,\,(z_T,z_{\Gamma,T})=(0,0)\,\,\;\textnormal{in}\,\,G\times\Gamma,\,\;\mathbb{P}\textnormal{-a.s.}
		\end{align*}
	\end{cor}
	\section{Controllability results}\label{sec4}
	This section establishes the null and approximate controllability results of the forward stochastic heat equation \eqref{1.1}. Let us first prove the null controllability result.
	\begin{proof}[Proof of Theorem \ref{thmm1.1}]
		Let $(y_{0},y_{\Gamma,{0}})\in L^2_{\mathcal{F}_{0}}(\Omega;\mathbb{L}^2)$ and consider the following linear subspace of $L_{\mathcal{F}}^{1}(0,T;L^2(\Omega;L^2(G)))$ defined by
		$$\mathcal{H}=\big\{\chi_E\chi_{G_0}z: \,\,\, (z,z_\Gamma,Z,\widehat{Z}) \,\,\text{solves}\,\, \eqref{backadj} \,\,\text{with}\,\, (z_{T}, z_{\Gamma, T})\in L^2_{\mathcal{F}_{T}}(\Omega;\mathbb{L}^2)\big\}.$$
		Let $h$ the linear functional defined on $\mathcal{H}$ by
		$$h(\chi_E\chi_{G_0}z)=-\mathbb{E}\int_G y_{0}\,z(s) \,dx-\mathbb{E}\int_\Gamma y_{\Gamma,0}\,z_{\Gamma}(s)\,d\sigma.$$
		By applying the observability inequality \eqref{1.301} with $s=0$, the linear functional $h$ is well defined, bounded on $\mathcal{H}$, and satisfies that
		$$\vert h\vert^2_{\mathcal{L}(\mathcal{H};\mathbb{R})}\leq C\vert(y_{0},y_{\Gamma,{0}})\vert^2_{L^2_{\mathcal{F}_{0}}(\Omega;\mathbb{L}^2)},$$
		where $C$ is the same observability constant as in \eqref{1.301}. Now, by Hahn-Banach Theorem, $h$ can be extended to a bounded linear functional $\tilde{h}$ on the space $L_{\mathcal{F}}^{1}(0,T;L^2(\Omega;L^2(G)))$ such that
		$\vert\tilde{h}\vert_{\mathcal{L}(L_{\mathcal{F}}^{1}(0,T;L^2(\Omega;L^2(G)));\mathbb{R})}=\vert h\vert_{\mathcal{L}(\mathcal{H};\mathbb{R})}$. Then, utilizing a Riesz-type Representation Theorem for stochastic processes (see, e.g., \cite[Theorem 2.73]{lu2021mathematical}), there exists a process $u\in L^{\infty}_\mathcal{F}(0,T;L^2(\Omega;L^2(G)))$ such that 
		\begin{align}\label{4.10141}
			\tilde{h}(\zeta)=\mathbb{E}\int_Q u \zeta\,dxdt,\quad\forall\zeta\in L_{\mathcal{F}}^{1}(0,T;L^2(\Omega;L^2(G))).
		\end{align}
		Moreover,
		\begin{align*}
			\vert u\vert^2_{L_{\mathcal{F}}^{\infty}(0,T;L^2(\Omega;L^2(G)))}\leq C\vert (y_{0},y_{\Gamma,{0}})\vert^2_{L^2_{\mathcal{F}_{0}}(\Omega;\mathbb{L}^2)}.
		\end{align*}
		We claim that the above-obtained process $u$ is our desired control. Now, from \eqref{4.10141}, we deduce particularly that for any $(z_{T}, z_{\Gamma, T})\in L^2_{\mathcal{F}_{T}}(\Omega;\mathbb{L}^2)$, we have
		\begin{equation}\label{Rez}
			\mathbb{E} \int_Q  \chi_E\chi_{G_{0}} u z\,dx \,dt+  \mathbb{E}\int_{G}y_{0}z(0) \,dx+\mathbb{E}\int_{\Gamma}y_{\Gamma,0}z_{\Gamma}(0) \, \,d\sigma =0. 
		\end{equation}
		From equations \eqref{1.1} and \eqref{backadj}, we obtain 
		\begin{equation}\label{Rez1}
			\mathbb{E}\int_{G}y(T)z_{T} \,dx+\mathbb{E}\int_{\Gamma}y_{\Gamma}(T)z_{\Gamma,T} \, \,d\sigma-\mathbb{E}\int_{G}y_{0}z(0) \,dx-\mathbb{E}\int_{\Gamma}y_{\Gamma,0}z_{\Gamma}(0) \, \,d\sigma= \mathbb{E} \int_Q  \chi_E\chi_{G_{0}} u z\,dx \,dt.
		\end{equation}
		Combining  \eqref{Rez} and \eqref{Rez1}, we conclude that
		$$\mathbb{E}\int_{G}y(T)z_{T} \,dx+\mathbb{E}\int_{\Gamma}y_{\Gamma}(T)z_{\Gamma,T} \, \,d\sigma = 0,$$
		which holds for all $(z_{T},z_{\Gamma,T})\in \mathbb{L}^2_{\mathcal{F}_{T}}(\Omega;\mathbb{L}^2)$. Therefore
		$$(y(T,\cdot),y_\Gamma(T,\cdot))=(0,0) \;\,\,\textnormal{in}\,\, G\times\Gamma,\,\,\,\mathbb{P}\textnormal{-a.s.}$$
		This completes the proof of Theorem \ref{thmm1.1}.
	\end{proof}
	Now, we address the proof of the approximate controllability property of \eqref{1.1}.
	\begin{proof}[Proof of Theorem \ref{thmm1.2}]
		Suppose that $\textbf{m}((s,T)\cap E)>0$, for any $s\in[0,T)$. Let us define the following operator 
		$$\mathcal{L}_T:L^\infty_\mathcal{F}(0,T;L^2(\Omega;L^2(G)))\longrightarrow L^2_{\mathcal{F}_T}(\Omega;\mathbb{L}^2),\;u\longmapsto (y(T,\cdot),y_\Gamma(T,\cdot)),$$
		where $(y,y_\Gamma)$ is the solution of \eqref{1.1} with control $u$ and initial condition $(y_0,y_{\Gamma,0})=(0,0)$. Then, it is easy to see that \eqref{1.1} is approximately controllable at time $T$ if and only if $\overline{\mathcal{R}(\mathcal{L}_T)}=L^2_{\mathcal{F}_T}(\Omega;\mathbb{L}^2)$ if and only if $\mathcal{L}^*_T$ is one-to-one. Let us now compute the adjoint operator $\mathcal{L}_T^*$. By Itô's formula, we calculate $d\langle(y,y_\Gamma),(z,z_\Gamma)\rangle_{\mathbb{L}^2}$ where $(y,y_\Gamma)$ is the solution of \eqref{1.1} and $(z,z_\Gamma,Z,\widehat{Z})$ is the solution of its adjoint equation \eqref{backadj} with $s=0$. Then, we obtain that
		\begin{align*}
			&\,\langle(y(T,\cdot),y_\Gamma(T,\cdot)),(z_T,z_{\Gamma,T})\rangle_{L^2_{\mathcal{F}_T}(\Omega;\mathbb{L}^2)}-\langle(y_0,y_{\Gamma,0}),(z(0,\cdot),z_\Gamma(0,\cdot))\rangle_{L^2_{\mathcal{F}_0}(\Omega;\mathbb{L}^2)}\\
			&=\mathbb{E}\int_Q \chi_E(t)\chi_{G_0}(x)uz \,dxdt\\
			&=\langle u,\chi_E \chi_{G_0}z\rangle_{L^2_\mathcal{F}(0,T;L^2(G))}\\
			&=\langle u,\chi_E \chi_{G_0}z\rangle_{L^\infty_\mathcal{F}(0,T;L^2(\Omega;L^2(G))),(L^\infty_\mathcal{F}(0,T;L^2(\Omega;L^2(G))))'}.
		\end{align*}
		It follows that $\mathcal{L}_T^*(z_T,z_{\Gamma,T})=\chi_E\chi_{G_0}z$. Hence, the approximate controllability of \eqref{1.1} is equivalent to the following property
		\begin{align*}
			z=0\,\,\;\textnormal{in}\,\,E\times G_0,\,\;\mathbb{P}\textnormal{-a.s.} \,\,\Longrightarrow\,\,\,(z_T,z_{\Gamma,T})=(0,0)\,\,\textnormal{in}\,\,G\times\Gamma\;\,\;\mathbb{P}\textnormal{-a.s.}
		\end{align*}
		From Corollary \ref{corunconprop}, the above uniqueness property holds, then the system \eqref{1.1} is approximately controllable.
		
		Conversely, suppose that \eqref{1.1} is approximately controllable.  By contradiction argument, assume that there exists $s_0\in[0,T)$ so that $\textbf{m}((s_0,T)\cap E)=0$. Let $\xi_1\in L^2_\mathcal{F}(s_0,T)$ be a nonzero process and set $(Z_1,\widehat{Z}_1)=\xi_1(\psi_1,\psi_{\Gamma,1})\in L^2_\mathcal{F}(s_0,T;\mathbb{L}^2)$ which is a nonzero process. Now, let us consider the following forward stochastic differential equation
		\begin{equation}\label{fordiff1}
			\begin{cases}
				\begin{array}{ll}
					d\varphi_1 - \lambda_1 \varphi_1 \,dt = -a(t)\varphi_1\,dt - b(t)\xi_1 \,dt + \xi_1 \,dW(t) &\textnormal{in}\,\,(s_0,T],\\
					\varphi_1(s_0)=0.
				\end{array}
			\end{cases}
		\end{equation}
		Noting that \eqref{fordiff1} has a nonzero unique solution $\varphi_1\in L^2_\mathcal{F}(\Omega;C([0,T];\mathbb{R}))$. Now, the following
		$$(z_1,z_{\Gamma,1},Z_1,\widehat{Z}_1)=(\varphi_1\psi_1,\varphi_1\psi_{\Gamma,1},\xi_1\psi_1,\xi_1\psi_{\Gamma,1})\in L^2_\mathcal{F}(\Omega;C([s_0,T];\mathbb{L}^2))\times L^2_\mathcal{F}(s_0,T;\mathbb{L}^2)$$
		is a nonzero solution to the following forward stochastic parabolic equation
		\begin{equation}\label{backadj3.10}
			\begin{cases}d z_1+\Delta z_1 d t=\big[-a(t) z_1-b(t) Z_1\big] d t+Z_1 d W(t) & \text { in }(s_0, T) \times G, \\ 
				d z_{\Gamma}+\Delta_{\Gamma} z_{\Gamma,1} d t-\partial_\nu z_1 d t=\big[-a(t) z_{\Gamma,1}-b(t) \widehat{Z}_1\big] d t+\widehat{Z}_1 d W(t) & \text { on }(s_0, T) \times \Gamma, \\ z_{\Gamma,1}(t, x)=\left.z_1\right|_{\Gamma}(t, x) & \text { on }(s_0, T) \times \Gamma, \\ \left.\left(z_1, z_{\Gamma,1}\right)\right|_{t=T}=\left(z_{1}^T,z_{\Gamma,1}^T\right) & \text { in } G \times \Gamma.
			\end{cases}
		\end{equation}
		Now, set 
		$$(z,z_\Gamma,Z,\widehat{Z})=\begin{cases}
			0 &\textnormal{in}\;\, (0,s_0),\\
			(z_1,z_{\Gamma,1},Z_1,\widehat{Z}_1) &\textnormal{in}\;\,(s_0,T).
		\end{cases}$$
		Hence, $(z,z_\Gamma,Z,\widehat{Z})$ is a nonzero solution of \eqref{backadj} for which $z(t,x)=0$ in $E\times G_0$, because of the assumption $\textbf{m}((s_0,T)\cap E)=0$. This contradicts the approximate controllability of \eqref{1.1}. Moreover, see also that $\left(z_1(s_0,\cdot), z_{\Gamma,1}(s_0,\cdot)\right)=(0,0)$, then in our case the backward uniqueness does not hold anymore for backward stochastic parabolic equations with dynamic boundary conditions. For this reason, the null controllability of \eqref{1.1} does NOT imply the approximate controllability as in the deterministic setting for heat equations. This concludes the proof of Theorem \ref{thmm1.2}.\\
	\end{proof}
	\section{Conclusion}\label{sec5}
	In this work, null and approximate controllability problems are addressed for the forward stochastic heat equation with dynamic boundary conditions. We first proved an appropriate observability inequality under a strong measurability condition by using a spectral inequality. Then, by applying the duality approach, our controllability results are established.
	
	Based on the results developed in this paper, many other interesting questions can be addressed. Here, we briefly describe some of them:
	\begin{itemize}
		\item [$\bullet$] When the coefficients $a$ and $b$ depend on the space variable, the controllability of \eqref{1.1} is an open problem. In general, we are interested in the study of the controllability of the following general forward stochastic parabolic system
		\begin{equation}\label{open1.1}
			\begin{cases}
				\begin{array}{ll}
					dy - \Delta y \,dt = \big[a_1y+\chi_{G_0}(x)u\big] \,dt + a_2y\,dW(t) &\textnormal{in}\,\,Q_0,\\
					dy_\Gamma-\Delta_\Gamma y_\Gamma \,dt+\partial_\nu y \,dt = b_1y_\Gamma \,dt+b_2y_\Gamma \,dW(t) &\textnormal{on}\,\,\Sigma_0,\\
					y_\Gamma(t,x)=y\vert_\Gamma(t,x) &\textnormal{on}\,\,\Sigma_0,\\
					(y,y_\Gamma)\vert_{t=0}=(y_0,y_{\Gamma,0}) &\textnormal{in}\,\,G\times\Gamma,
				\end{array}
			\end{cases}
		\end{equation}
		with $(y_0,y_{\Gamma,0})\in L^2_{\mathcal{F}_0}(\Omega;\mathbb{L}^2)$ is the initial state, $u\in L^2_{\mathcal{F}}(0,T;L^2(G))$ is the control function, $a_1,a_2\in L^\infty_{\mathcal{F}}(0,T;L^\infty(G))$ and $b_1,b_2\in L^\infty_{\mathcal{F}}(0,T;L^\infty(\Gamma))$. The main difficulty in establishing the null controllability of \eqref{open1.1} is to show the following challenging observability inequality: There exists a constant $C>0$ so that
		$$\mathbb{E}|z(0)|^2_{L^2(G)}+\mathbb{E}|z_\Gamma(0)|^2_{L^2(\Gamma)}\leq C\,\mathbb{E}\int_0^T\int_{G_0}z^2 \,dxdt,$$
		for all $(z,z_\Gamma,Z,\widehat{Z})$ solution of the corresponding adjoint backward stochastic parabolic system. One way to solve it is to construct a suitable weight function allowing us to remove the correction terms ``$Z$'' and ``$\widehat{Z}$'' on the right-hand side in the  Carleman estimate. This problem is still open for Dirichlet boundary conditions as well, see, e.g., \cite[Chapter 9]{lu2021mathematical} and \cite{tang2009null} for more details.
		\item [$\bullet$] The controllability of \eqref{open1.1} without the diffusion on the boundary is also an interesting question. In \cite{elgrouDBC}, such a term was essential to derive the appropriate Carleman estimate for backward stochastic parabolic equations with dynamic boundary conditions. It appears that the controllability of such a problem in the deterministic setting is an open problem as well.
		\item [$\bullet$] The controllability of semilinear and quasi-linear stochastic parabolic equations is far from being studied. To the author's best knowledge, there are no results in the quasi-linear setting. For semi-linear systems with global Lipshitz non-linearity, we refer to \cite{san23}. 
		\item  [$\bullet$] It seems that the observability inequality \eqref{1.301} will play a crucial role in the study of time optimal control issues for the forward stochastic parabolic equation \eqref{1.1}. This will be the subject of our forthcoming works. For some details on the time and norm optimal control problems for (deterministic) parabolic equations with dynamic boundary conditions, we refer the readers to \cite{OuBoMa21}.
	\end{itemize}

\end{document}